\documentclass[11pt]{article}
\author{Bing-Long Chen\\[8pt]}
\title{\textbf{On stationary solutions to the non-vacuum  Einstein field equations}}
\date{June 14, 2016}
\usepackage{latexsym}
\usepackage{amsmath}
\usepackage{amssymb,amscd}
\usepackage[mathscr]{eucal}
\newtheorem{thm}{Theorem}[section]

\newtheorem{lem}[thm]{Lemma}

\newtheorem{rem}{Remark}[section]

\numberwithin{equation}{section}
\newenvironment{pf}{{\noindent \it  Proof.}}{{\hfill$\Box$}\\}
\begin{document}
\maketitle
\let\thefootnote\relax\footnotetext{AMS Mathematics Subject Classification Numbers: Primary
53c50; Secondary 83c20. }

\begin{abstract}
We derive a local curvature estimate for four-dimensional  stationary solutions to the inheriting Einstein-Maxwell-Klein-Gordon  equations.  In particular, it implies that any such stationary geodesically complete solution with vanishing Poynting vector and  proper coupling constants (like dark energy) is flat.  We also generalize the results to  higher dimensions. 
\end{abstract}

  \section{Introduction} 
  The purpose  of this paper  is to generalize the results in \cite{Chen1} for vacuum Einstein field equations  to non-vacuum Einstein field equations. In \cite{Chen1}, we proved that any geodesically complete 4-d spacetime with a timelike Killing field satisfying the vacuum Einstein field equation with nonnegative cosmological constant is  actually flat. This result generalized the previous  results in \cite{A} and \cite{L} by removing the extra conditions (like chronological condition)  imposed on their theorems.   
  
  Because the physical spacetime can never be flat, it must violate the assumptions in the theorem in \cite{Chen1}, for instance, the geodesic completeness.  Also, by the singularity theorem of Penrose and Hawking (see \cite{HE} Chapter 8), the geodesic completeness of the physical spacetime is also hardly valid.  Moreover,  we know that many stationary spacetimes contain black holes, see \cite{CLH} \cite{Ro}.    
  In \cite{Chen1},  in order to prove the above result, we actually derived a local curvature estimate (see Theorem 1.3 in \cite{Chen1}), which does  not require the spacetime to be geodesically complete. So from this point of view, the local curvature estimate in \cite{Chen1}  seems to be more important.

 In this paper,  we are mainly concerned with the Einstein field equations  in the presence of two kinds of matters or fields,  the electromagnetic fields and scalar fields. For the former case, the interaction of the gravity (a 4-d spacetime $(M,g_{M})$)  and the electromagnetic field  (a 2-form $F$)  is described by a system of differential equations  called Einstein-Maxwell equation: 
  \begin{equation} \label{EM}
  \begin{split}
  R_{\alpha\beta}-\frac{1}{2}Rg_{\alpha\beta}+\Lambda g_{\alpha\beta}&=\kappa (F_{\alpha\gamma}F_{\beta\delta}g^{\gamma\delta}-\frac{1}{4}|F|^2 g_{\alpha\beta})\\
   dF&=d\ast F=0,
  \end{split}
  \end{equation}
     where  $\Lambda, \kappa$ are constants, $|F|^2\triangleq F_{\xi\gamma}F_{\eta\delta}g^{\xi\eta}g^{\gamma\delta}$. The second equation in (\ref{EM}) satisfied by the 2-form $F_{\alpha\beta}$ (electromagnetic field) is  called the Maxwell  equation.

A scalar field, in physics,  associates a value  ($\in R$  or $C$) to each point of the spacetime, like temperature, pressure etc.   In this  paper,  we consider the  Klein-Gordon scalar fields,  which can be used to describe the $\pi^{0}$, or $\pi^{+},\pi^{-}$ mesons.  The Einstein equation coupled with Klein-Gordon equation is the following 
 \begin{equation}\label{EKG}
 \begin{split}
&Ric(g_{M})-\frac{R}{2}g_M+\Lambda g_M=\kappa'[d\phi\otimes d\phi-\frac{1}{2}(|d\phi|^2+\frac{m^2}{\hbar^2}|\phi|^2)g_{M}]\\
&\triangle_{g_M} \phi =\frac{m^2}{\hbar^2} \phi
\end{split}
\end{equation}
where $\Lambda,\kappa',$ $m$ (mass), $\hbar$ (reduced  planck constant) are  constants, $\phi: M\rightarrow R$ or $R^2$ is a map (scalar field).

The (uncharged) Einstein-Maxwell-Klein-Gordon equation is a combination of (\ref{EM}) and (\ref{EKG}) (see \cite{HE} Chapter 3): 
 \begin{equation}\label{EMKG}
 \begin{split}
&Ric(g_{M})-\frac{R}{2}g_M+\Lambda g_M=\kappa(F_{\alpha\gamma}F_{\beta\delta}g^{\gamma\delta}-\frac{1}{4}|F|^2g_{\alpha\beta})\\& \ \ \ \ \ \ \ \ \ \ \  \ \ \ \ \ +\kappa'[d\phi\otimes d\phi-\frac{1}{2}(|d\phi|^2+\frac{m^2}{\hbar^2}|\phi|^2)g_{M}]\\
&dF=d\ast F=0\\
&\triangle_{g_M} \phi =\frac{m^2}{\hbar^2} \phi.
\end{split}
\end{equation}
 
It is easy to see that (\ref{EMKG}) reduces to (\ref{EM}) or (\ref{EKG}) when $\phi=0$ or $F=0$ in (\ref{EMKG}). 
 So instead of dealing with (\ref{EM}) and (\ref{EKG}) separately, we can handle only one equation (\ref{EMKG}).

    Let $X$ be a timelike Killing field on $M$,  by which we say the spacetime $(M,g_{M})$ is stationary.   Following \cite{To1} \cite{To2},  we say the solution $(g_{M}, F,\phi)$ in (\ref{EMKG}) is inheriting if $F$ and $\phi$ satisfy   \begin{equation} \label{1.4}
  \mathcal{L}_X(F)=0, \ d\phi(X)=0, 
 \end{equation}
which means that  $F$ and $\phi$ inherit the $X$-symmetry  from the metric tensor $g_{M}$.

 Denote the electric and  magnetic fields (related to $X$) by \begin{equation}\label{EB}E=i_{X}F, B=i_{X}\ast F\end{equation} respectively, where $i_{X}F$ is the tensor contraction of $F$ with  $X$.   Let $X^{\ast}$ be the  1-form on $M$ obtained from $X$ by lowering indices.   We associate   a Riemannian metric tensor $\hat{g}$ on $M$ to $g_{M}$ (see \cite{Chen1}):  \begin{equation}\label{hat}
\hat{g}\triangleq -\frac{2}{g_{M}(X,X)}X^{\ast} \otimes X^{\ast}+g_{M}.
\end{equation}

The first main result of the paper is  the following local curvature estimate:  
  \begin{thm} \label{N1.1}  Let $(M,g_{M})$ be a 4-d spacetime with a timelike Killing field $X$ such that $g_{M}$ satisfies the inheriting Einstein-Maxwell-Klein-Gordon  equation (\ref{EMKG}) (\ref{1.4}),  where $\Lambda\geq 0$ and $\kappa\leq 0,\kappa'\geq  0$ are constants.    Let $\hat{B}(x_0, a)$ be a $\hat{g}-$metric ball centered at $x_0$ of radius $a>0$ with compact  closure  in $M$ and  \begin{equation} \label{1.7}
\begin{split}
\sup_{x\in \hat{B}(x_0,{a})}|\kappa|u^{-2}\mid E\wedge B \mid_{\hat{g}}\leq a^{-2}.
\end{split}
\end{equation}
 Then for any $0<\delta<1$, there is a constant $C_{\delta}>0$ depending only on $\delta$  such that 
 \begin{equation} \label{1.8}
\begin{split}
\sup_{x\in \hat{B}(x_0,\frac{a}{2})}\mid Rm(g_M)\mid_{\hat{g}}\leq \delta \frac{m^2}{\hbar^2}+C_{\delta}a^{-2}.
\end{split}
\end{equation}
\end{thm}

    The Hodge dual (on three-space) of the 2-form  $E\wedge B$ in (\ref{1.7}) is a  1-form corresponding to the Poynting vector in physics.         
     A corollary of Theorem \ref{N1.1} is the following 
      \begin{thm} \label{N1.2}  Let $(M,g_{M})$ be a geodesically complete 4-d spacetime with a timelike Killing field $X$ such that $g_{M}$ satisfies the inheriting Einstein-Maxwell-Klein-Gordon  equation (\ref{EMKG}) (\ref{1.4}),  where $\Lambda\geq 0, \kappa\leq 0, \kappa'\geq 0$.    We assume the Poynting vector vanishes on $M$, i.e., $E\wedge B\equiv 0$. Then $(M,g_{M})$  is flat. 
       \end{thm}
       We give a  possible physical explanation why Theorem \ref{N1.1} or Theorem \ref{N1.2} hold  for  $\Lambda\geq 0, \kappa\leq 0, \kappa'\geq 0$ (but not the contrary). It is well-known that positive cosmological constant $\Lambda$ and negative energy ($\kappa<0$, $\kappa'>0$)  will cause  the acceleration of the universe. Intuitively, the acceleration  will reduce the density of the matter or energy distribution of the universe,   hence make the spacetime more "regular".  
       
       \begin{rem} The spacetime being static   implies  that  the Poynting vector vanishes, i.e. $E \wedge B=0$.  
\end{rem}
We say a stationary spacetime is static if the orthogonal complement  of the timelike Killing field is an integrable distribution. In general, when dim $\geq 5$, the "curvature" measuring the difference between stationaryness  and staticity is difficult to control. In order to generalize Theorems  \ref{N1.1} and \ref{N1.2} to dim $\geq 5$, we can  assume for simplicity that the spacetime is static. In this case,  we can obtain  analogous results for general dimensions,  see Theorems \ref{N3.1}, \ref{N3.2}.

 The paper is organized as follows. In section 2,   we prove Theorems \ref{N1.1} and \ref{N1.2}. In section 3, we generalize the results to general  dimensions.  
 
 \textbf{Acknowledgement } The author is grateful to Professors  S. T. Yau, X. P. Zhu and Dr. J. B. Li for helpful discussions.   
 The work was partially supported  by grants NSFC11521101, 11025107.

 \section{4-d Einstein-Maxwell-Klein-Gordon  equations} 
 \subsection{ Preliminaries}
  In this section, we recall some useful formulas for stationary spacetimes, their  computations could be found in \cite{Chen1},  see also \cite{A} \cite{KSMH}. 
  Let $\{x^{\alpha}\}$ be a local coordinate system on a 4-d  stationary spacetime $(M,\bar{g}, X)$ so that the metric tensor $\bar{g}=g_{\alpha\beta}dx^{\alpha}dx^{\beta}$ takes the following form: 
     \begin{equation}\label{2.1}
     \bar{g}=-u^{2}(dt+\theta)^2+g
     \end{equation}
     where $u=(-\langle X,X\rangle)^{\frac{1}{2}}$, $\theta=\theta_idx^i$ and $g=g_{ij}dx^idx^j$ are $t=x^0$ -independent. Here we use $\bar{g}$ to denote the spacetime metic $g_{M}$. We use the convention that for tensor indices,  the Greek letters $\alpha, \beta,\cdots$ are running from $0$ to $3$, the Latin letters $i,j,\cdots$ from $1$ to $3$.  Let $e_0=\frac{\partial}{\partial t}$,  $ e_i= \frac{\partial}{\partial x^i}-\theta_i \frac{\partial}{\partial t}$ be a local tangent frame satisfying $\langle e_i,e_j\rangle=g_{ij}$ and $\langle e_0,e_j\rangle=0$. Fix the  local orientation of $M$ (on this coordinate chart) given by the frame $(e_1,e_2,e_3,e_0)$. We use $\ast^4$ and $\ast$ to denote the Hodge dual operators (with respect to the above orientation and metics $\bar{g}$ and $g$) on spaces of  dimensions  4 and 3 respectively. 
     
     We remark that the coordinate system $\{x^{\alpha}\}$ and the local tangent frame $\{e_{\alpha}\}$ will be  used frequently through out of the paper. It is clear that their  constructions are  irrelevant to  the dimensions. 
     
     Let $\omega=u^{3}\ast d\theta$ be a spatial 1-form on $M$.   For any  Lorentzian metric $\bar{g}$ of the form (\ref{2.1}),  the Ricci curvature components $\bar{R}ic(e_{\alpha},e_{\beta})$ can be expressed as a system of differential equations satisfied by $u, \omega$ and $g$ (see \cite{Chen1}, (5.1)): 
     \begin{equation} \label{2.2}
\begin{split}
& R_{ij}={u^{-1}}{\nabla_{i}\nabla_j u}+\frac{1}{2}{u^{-4}}(\omega_i\omega_j-|\omega|^2 g_{ij})+\bar{R}ic(e_i,e_j)\\
& \triangle u =-\frac{1}{2}u^{-3}|\omega|^2+ u^{-1}\bar{R}ic(X,X) \\
& g^{kl} \nabla_k\omega_{l}=3g^{kl}\omega_{k}\nabla_{l}\log u\\
& (\ast d\omega)_j=\pm2u \bar{R}ic(X,e_j).
\end{split}
\end{equation} 
Let $\tilde{g}\triangleq u^2g$ be a conformal change of the spatial metric $g$,  the Ricci curvature of $\tilde{g}$  can be computed(see \cite{Chen1}, (5.3)):   \begin{equation} \label{2.3}
\begin{split}
& \tilde{R}_{ij}=\frac{1}{2}u^{-4}\omega_{i}\omega_{j}+2\frac{u_iu_j}{u^2}+\bar{R}ic(e_i,e_j)-{u^{-2}\bar{R}ic(X,X)}g_{ij}.
\end{split}
\end{equation}
The Ricci curvatures  of the metrics $\hat{g}$ (see (\ref{hat})) and $\bar{g}$ are related (see \cite{Chen1}, (5.2)) in the following manner: 
 \begin{equation} \label{2.4}
\begin{split}
& \hat{R}ic({X,X})=u^{-2}| \omega|^2-\bar{R}ic(X,X) \\
& \hat{R}ic({X,e_j})=-\bar{R}ic(X,e_j)
\\& \hat{R}ic({e_i,e_j})=-u^{-4}(|\omega|^2g_{ij}-\omega_i\omega_j)+\bar{R}ic(e_i,e_j).
\end{split}
\end{equation} 
For later use,   the $\hat{g}$-hessian of 
any time-independent smooth function $f$ on the spacetime can be computed to be (see \cite{Chen1}, (2.11)): 
 \begin{equation} \label{2.5}
\begin{split}
& \hat{\nabla}^2f(e_0,e_0)=\frac{1}{2}\langle\nabla u^2,\nabla f\rangle \\
&\hat{\nabla}^2f(e_0,e_j)=\frac{1}{2}u^{-1} (\ast (\omega\wedge df))(e_j)\\
& \hat{\nabla}^2f(e_i,e_j)=\nabla_{ij}f.
\end{split}
\end{equation} In particular, we have        
 \begin{equation} \label{2.6}
\begin{split}
& \hat{\triangle}f=\bar{\triangle}f=u^2\tilde{\triangle}f= \triangle f+\langle \nabla \log u,\nabla f\rangle.
\end{split}
\end{equation}
 \subsection{Field equations}
 
 Denote $V\circ \phi = \frac{m^2}{\hbar^2}|\phi|^2$.  Taking trace on  (\ref{EMKG}), we find $$R=\kappa' |d\phi|^2+4\Lambda+2\kappa' V$$ and equation (\ref{EMKG} ) can be rewritten as 
\begin{equation}\label{2.7}
R_{\alpha\beta}=\kappa (F_{\alpha\gamma}F_{\beta\delta}g^{\gamma\delta}-\frac{1}{4}|F|^2 g_{\alpha\beta})+\kappa' \phi_{\alpha}\phi_{\beta}+(\Lambda+\frac{\kappa'}{2}V)g_{\alpha\beta}.
\end{equation}

  From $\mathcal{L}_{X}F=0$,  one can derive  $\mathcal{L}_{X}\ast^4  F=0$ since $X$ is a Killing field. Denote the electric and  magnetic fields related to $X$  by $$E=i_{X}F, B=i_{X}\ast^4 F$$ respectively.  From Cartan's  homotopy formula $\mathcal{L}_{X}=di_{X}+i_{X}d$ for differential forms, we find $dE=dB=0$.  Using  $i^2_{X}=0$, we obtain \begin{equation}\label{2.8}
     \mathcal{L}_{X}E=\mathcal{L}_{X}B=0, i_{X}(E)=i_{X}(B)=0,
     \end{equation} which roughly say that $E$ and $B$ are spatial 1-forms on $M$.

      Let $F=\omega_{1}(x)\wedge u(dt+\theta)+\omega_2(x)$ be a local expression of the electromagnetic field $F$ in the coordinate system $\{x^{\alpha}\}$ in Section 2.1,  where $\omega_{1}(x)$ and $\omega_2(x)$ are spatial forms of degrees  1 and 2. Since  $\ast^4F=\ast \omega_2\wedge u(dt+\theta)+\ast \omega_1$, we know $ E=-u\omega_1(x)$, $B(x)=-u\ast \omega_2$.  From this, we  know $F_{\alpha\beta}F_{\gamma\delta}g^{\alpha\gamma}g^{\beta\delta}=2u^{-2}(-|E|^2+|B|^2)$. The Maxwell's  equation $dF=d\ast^4F=0$ is   equivalent to  a system of differential equations on the electric and magnetic fields $E$ and $B$:
        \begin{equation} \label{2.9}
  \begin{split}
  & dE=0\\
  & dB=0\\
 & d(u^{-1}\ast E)= B \wedge d\theta\\
 & d (u^{-1}\ast  B)=E\wedge d\theta.
  \end{split}
  \end{equation}
 Now we compute the components of the Ricci curvature $\bar{R}ic$ in terms of the frame $\{e_{\alpha}\}$ in Section 2.1.  From (\ref{2.7}),  we have 
 \begin{equation}\label{2.10}
 \begin{split} \bar{R}ic(X,X)&=\kappa(\frac{1}{2}|E|^2+\frac{1}{2}|B|^2)-(\Lambda+\frac{\kappa'}{2}V)u^2.
  \end{split}
 \end{equation}
     Since $i_{e_j}F=-E(e_j)(dt+\theta)+i_{e_j}\omega_2$ and $$\langle i_{e_i}\omega_2, i_{e_j}\omega_2\rangle=|\ast \omega_2|^2g_{ij}-(\ast \omega_2)(e_i)(\ast \omega_2)(e_j),$$
     we have 
     $$
     \langle i_{e_i}F, i_{e_j}F\rangle=-u^{-2}(E_iE_j+B_iB_j)+u^{-2}|B|^2 g_{ij}.
     $$
     
     This implies 
     \begin{equation} \label{2.11}
     \begin{split}
 \bar{R}ic(e_i,e_j)
 =&\kappa[-u^{-2}(E_iE_j+B_iB_j)+\frac{1}{2}u^{-2}(|B|^2+|E|^2)g_{ij}]\\&+\kappa'\phi_i\phi_j+(\Lambda+\frac{\kappa'}{2}V)  g_{ij}.
 \end{split}
  \end{equation}

  Hence  \begin{equation}\label{2.12}
\begin{split}
&\bar{R}ic(e_k,e_l)-u^{-2}\bar{R}ic(X,X)g_{kl}\\& = -\kappa u^{-2}(E_iE_j+B_iB_j)+\kappa' \phi _k\phi _l+(2\Lambda+{\kappa'}V)g_{kl}.
\end{split}
\end{equation}

    To compute  the term  $\bar{R}ic(X,e_j)$,  we need 
  \begin{equation}\label{2.13}
  \begin{split}
     \langle i_{X}F, i_{e_j}F\rangle&=\langle E, -E(e_j)(dt+\theta)+i_{e_j}(\omega_2) 
     \rangle\\
     &=\langle E, i_{e_j}(\omega_2) \rangle\\
     &=i_{e_{j}}\ast (E\wedge \ast \omega_2)\\
     &=-u^{-1} i_{e_j} (\ast (E\wedge B)).
     \end{split}
     \end{equation}
     Therefore, 
     \begin{equation}\label{2.14}
  \begin{split}
     \bar{R}ic(X,e_j)=-\kappa u^{-1} i_{e_j} (\ast (E\wedge B)). 
     \end{split}
     \end{equation}
  Combining (\ref{2.2}) (\ref{2.3}) (\ref{2.11}) and  (\ref{2.12}),  we have 
 \begin{equation} \label{2.15}
\begin{split}
 R_{ij}=& {u^{-1}}{\nabla_{i}\nabla_j u}+\frac{1}{2}{u^{-4}}(\omega_i\omega_j-|\omega|^2 g_{ij})+\kappa'\phi_i\phi_j+(\Lambda+\frac{\kappa'}{2}V) g_{ij}\\&+\kappa[-u^{-2}(E_iE_j+B_iB_j)+\frac{1}{2}u^{-2}(|B|^2+|E|^2)g_{ij}],
\end{split}
\end{equation} 
and 
\begin{equation} \label{2.16}
\begin{split}\tilde{R}_{ij}=& 2u^{-2}u_iu_j+\frac{1}{2}u^{-4}\omega_i\omega_j-\kappa u^{-2}(E_iE_j+B_iB_j)\\
& +\kappa' \phi _i\phi _j+(2\Lambda+{\kappa'}V)g_{ij},
\end{split}
\end{equation}
  where $\tilde{R}_{ij}$ is the Ricci curvature of  $\tilde{g}=u^2g$.

 Denote 
 \begin{equation}\label{2.17}
 \begin{split}
 E=d\phi_3, \ \ \ B=d\phi_4
 \end{split}
 \end{equation}
  for some locally defined functions $\phi_3$ and $\phi_4$.  From (\ref{2.2}) (\ref{2.9}) (\ref{2.10}),  we have: 
     \begin{equation}\begin{split}\label{2.18}
      \triangle u&=-\frac{1}{2}u^{-3}|\omega|^2+\kappa (\frac{1}{2} u^{-1}|E|^2+\frac{1}{2}u^{-1}|B|^2)-(\Lambda+\frac{\kappa'}{2}V) u\\
     \triangle \phi_3&=\langle d\log u, E\rangle+u^{-2}\langle \omega,B\rangle\\
      \triangle \phi_4&=\langle d\log u, B\rangle+u^{-2}\langle \omega,E\rangle.
     \end{split}
      \end{equation}
    Combining (\ref{2.6}) and (\ref{2.18}),  we have 
 \begin{equation}\label{2.19}
     \begin{split}
      \hat{\triangle} \log u&=u^2\tilde{\triangle}\log u=-\frac{1}{2}u^{-4}|\omega|^2+\kappa (\frac{1}{2} u^{-2}|E|^2+\frac{1}{2}u^{-2}|B|^2)-(\Lambda+\frac{\kappa'}{2}V)\\
    \hat{\triangle} \phi_3&=u^2\tilde{\triangle}\phi_3=2\langle d\log u, E\rangle+u^{-2}\langle \omega,B\rangle\\
     \hat{\triangle} \phi_4&=u^2\tilde{\triangle}\phi_4=2\langle d\log u, B\rangle+u^{-2}\langle \omega,E\rangle.
     \end{split}
      \end{equation}

\subsection{A Bochner formula}
In  this section, we assume that the condition \begin{equation}\label{2.20}\bar{R}ic(X,e_j)=0
\end{equation}  always holds.  When condition (\ref{2.20})  holds, from the last 2 equations of (\ref{2.2}), we know  $\omega=d\psi$  holds locally for some function $\psi$   and \begin{equation}\label{2.21}
\hat{\triangle}\psi=u^2\tilde{\triangle} \psi=4\langle \omega,d\log u\rangle.
\end{equation}
  From (\ref{2.14}),  the   condition (\ref{2.20}) is equivalent to \begin{equation}\label{2.22}E\wedge B=0,  \end{equation} i.e., the Poynting vector vanishes.   Here we assume $E=B=0$  if $\kappa=0$. 
    
  \begin{lem} \label{l2.1}Let $g_{H}=y_1^{-2}dy_1^2+c_2y_1^{-l_2}dy_2^2+\cdots+c_my_1^{-l_m}dy^2_m$, $c_a>0$ be a Riemannian metric on m-dimensional upper half space $H=\{(y_1,y_2,\cdots,y_m): y_i\in R, y_1>0\}$.  In terms of  the natural frame $\{\frac{\partial}{\partial y^a}\}$, the only nonzero components of the Christoffel symbols  $\Gamma^{a}_{bc}$ of $g_{H}$ are: 
\begin{equation} \label{2.23}
\begin{split}
& {\Gamma}^{1}_{bb}=\frac{1}{2}c_bl_by_1^{1-l_b},   {\Gamma}^{b}_{1b}=-\frac{1}{2}l_by_1^{-1},  {\Gamma}^{1}_{11}=-y_1^{-1}, 
\end{split}
\end{equation}
where $b\neq 1$. Up to the symmetries of curvature tensors, the only nonzero components of the curvature tensor of $g_{H}$ are 
\begin{equation} \label{2.24}
\begin{split}
& R_{1a1a}=-\frac{1}{4}c_al_a^2y_1^{-l_a-2}, R_{abab}=-\frac{1}{4}c_ac_bl_al_by_1^{-l_a-l_b},
\end{split}
\end{equation}
where $a,b\neq 1$ and $a\neq b$. 
\end{lem}
   
   In particular, the sectional curvatures of  the metric $g_{H}$ in Lemma \ref{l2.1} are  nonpositive. 

Let $g_{H}=y_1^{-2}{(dy_1^2+dy_2^2-2\kappa y_1dy_3^2-2\kappa y_1 dy_4^2)}$ be a  metric on 4-d upper half space $H=\{(y_1,y_2,y_3,y_4): y_i\in R, y_1>0\}$.  We define a map $\varPhi$ from the coordinate system $\{x^{\alpha}\}$ of Section 2.1 to $H$ to be  $\varPhi=(u^2,\psi,\phi_3,\phi_4)=(y_1,y_2,y_3,y_4)$.  Because $u,\psi,\phi_3,\phi_4$  are time-independent, $\varPhi$ is also time-independent, hence  a map from $\{x^{i}\}$ to $H$. 

\begin{lem} \label{l2.2}The map  $\varPhi$  satisfies \\ 
\ \ \ \ \ \ i)  $\varPhi^{\ast} g_{H}=4d\log u\otimes d\log u+u^{-4}\omega\otimes\omega-2\kappa u^{-2} E\otimes E-2\kappa u^{-2}B\otimes B$; \\
\ \ \ \ \ \ ii)  $\hat{\triangle}\varPhi=u^{2}\tilde{\triangle}\varPhi=-(2\Lambda+\kappa'V) u^2 \frac{\partial}{\partial y_1}+u^{-2}\langle B,\omega\rangle\frac{\partial}{\partial y_3}+u^{-2}\langle E,\omega\rangle\frac{\partial}{\partial y_4}$, where $\hat{\triangle}\varPhi$ (or $\tilde{\triangle}\varPhi$) is the  harmonic map Laplacian of $\varPhi$ from  Riemannian manifold $(\{x^{\alpha}\}, \hat{g})$(or $(\{x^i\}, \tilde{g})$) to  $(H,g_{H})$.
\end{lem}
\begin{pf} Substituting $m=4$, $c_2=1,c_3=-2\kappa,c_4=-2\kappa$ and $l_{2}=2,l_3=1,l_4=1$ in Lemma \ref{l2.1}, we find
\begin{equation} \label{2.25}
\begin{split}
{\Gamma}^{1}_{11}=-y_1^{-1},  {\Gamma}^{1}_{22}=y_1^{-1}, {\Gamma}^{1}_{33}= {\Gamma}^{1}_{44}=-\kappa,  {\Gamma}^{2}_{12}=-y_1^{-1}, {\Gamma}^{3}_{13}=\Gamma_{14}^{4}=-\frac{1}{2}y_1^{-1}.
\end{split}
\end{equation}
By direct computations, 
 \begin{equation} \label{2.26}
\begin{split}
 \Gamma^{1}_{ab} \nabla_{\alpha}\varPhi^{a}\nabla_{\beta}\varPhi^{b}=&\Gamma^{1}_{11}\partial_{\alpha} u^2 \partial_{\beta} u^2+\sum_{b\neq 1}\Gamma^{1}_{bb} \partial_{\alpha} y_b \partial_{\beta} y_b\\
 =&-4u_{\alpha}u_{\beta}+u^{-2}\omega_{\alpha}\omega_{\beta}-\kappa (\phi_3)_{\alpha}(\phi_3)_{\beta}-\kappa (\phi_4)_{\alpha}(\phi_4)_{\beta}\\ 
  \Gamma^{2}_{ab} \nabla_{\alpha}\varPhi^{a}\nabla_{\beta}\varPhi^{b}=&   \Gamma^{2}_{12}( \nabla_{\alpha}\varPhi^{1}\nabla_{\beta}\varPhi^{2}+ \nabla_{\alpha}\varPhi^{2}\nabla_{\beta}\varPhi^{1})\\
  =& -2\omega_{\beta} (\log u)_{\alpha} -2\omega_{\alpha} (\log u)_{\beta}\\
   \Gamma^{3}_{ab} \nabla_{\alpha}\Phi^{a}\nabla_{\beta}\varPhi^{b}=&   \Gamma^{3}_{13}( \nabla_{\alpha}\varPhi^{1}\nabla_{\beta}\varPhi^{3}+ \nabla_{\alpha}\varPhi^{3}\nabla_{\beta}\varPhi^{1})\\
  =& -(\phi_3)_{\beta} (\log u)_{\alpha} -(\phi_3)_{\alpha} (\log u)_{\beta}\\
    \Gamma^{4}_{ab} \nabla_{\alpha}\varPhi^{a}\nabla_{\beta}\varPhi^{b}  =& -(\phi_4)_{\beta} (\log u)_{\alpha} -(\phi_4)_{\alpha} (\log u)_{\beta}.
\end{split}
\end{equation} 
Using (\ref{2.26}) and the formula   
 \begin{equation} \label{2.27}
\begin{split}
& (\hat{\nabla}_{\alpha\beta} \varPhi)^{a}= \hat{\nabla}_{\alpha\beta} \varPhi ^{a}+\Gamma^{a}_{bc} \nabla_{\alpha}\varPhi^{b}\nabla_{\beta}\varPhi^{c},
\end{split}
\end{equation}
we obtain
\begin{equation} \label{2.28}
\begin{split}
& (\hat{\nabla}_{\alpha\beta} \varPhi)^{1}= 2u^2 \hat{\nabla}_{\alpha\beta}\log u+u^{-2}\omega_{\alpha}\omega_{\beta}-\kappa E_{\alpha}E_{\beta}-\kappa B_{\alpha}B_{\beta}\\
& (\hat{\nabla}_{\alpha\beta} \varPhi)^{2}=  \hat{\nabla}_{\alpha\beta} \psi-2\omega_{\alpha}(\log u)_{\beta}-2\omega_{\beta}(\log u)_{\alpha}\\
& (\hat{\nabla}_{\alpha\beta} \varPhi)^{3}=\hat{\nabla}_{\alpha\beta} \phi_3-(\log u)_{\alpha}E_{\beta}-(\log u)_{\beta}E_{\alpha}\\
& (\hat{\nabla}_{\alpha\beta} \varPhi)^{4}=\hat{\nabla}_{\alpha\beta} \phi_4-(\log u)_{\alpha}B_{\beta}-(\log u)_{\beta}B_{\alpha}.
\end{split}
\end{equation}

Taking traces  on (\ref{2.28}) with respect to $\hat{g}$,  we get
\begin{equation}\label{2.29}
\begin{split}
& (\hat{\triangle}\varPhi)^{1}= 2u^2 \hat{\triangle}\log u+ u^{-2}|\omega|^2_{{g}}-\kappa|E|^2-\kappa|B|^2=-(2\Lambda+\kappa'V) u^2\\
& (\hat{\triangle} \varPhi)^{2}=  \hat{\triangle}\psi-4\langle \nabla \psi, \nabla \log u\rangle_{{g}}=0\\
& (\hat{\triangle} \varPhi)^{3}=\hat{\triangle}\phi_3-2\langle d\log u, E \rangle_{g}=\langle u^{-2}\omega, B\rangle_{g}\\
& (\hat{\triangle} \varPhi)^{4}=\hat{\triangle}\phi_4-2\langle d\log u, B \rangle_{g}=\langle u^{-2}\omega, E\rangle_{g},
\end{split}
\end{equation}
where we have used (\ref{2.19}) (\ref{2.21}). 
\end{pf}

\begin{lem}
Under the condition (\ref{2.20}) or (\ref{2.22}), we have $\langle\omega,E\rangle=\langle\omega, B\rangle=0$. Hence, $\hat{\triangle}\varPhi=-(2\Lambda+\kappa'V) u^2\frac{\partial}{\partial y_1}$.  
\end{lem}
\begin{pf} 
Assume $\kappa\neq 0$, otherwise the result holds trivially because we can take $E=B=0$. In order to show $\langle\omega,E\rangle=\langle\omega, B\rangle=0$, we argue by contradiction. Suppose $\langle\omega,E\rangle\neq 0$ at some point $P$,   by (\ref{2.22}), we have $B=\lambda E$ holds for some function $\lambda$ near $P$. Using i) in Lemma \ref{l2.2} and the contracted 2nd Bianchi identity, 
\begin{equation}\label{2.30}
\begin{split}
\frac{1}{2}\tilde{R}_i &=\tilde{\nabla}_j\tilde{R}_{ik}\tilde{g}^{kj}=\tilde{\nabla}_j(\frac{1}{2}(\varPhi^{\ast}h)_{ik}+\kappa'\phi_i\phi_k+(2\Lambda+\kappa'V) g_{ik})\tilde{g}^{jk}\\
&=\frac{1}{4} (tr_{\tilde{g}}(\varPhi^{\ast}h))_i +\frac{1}{2}\langle\tilde{\triangle}\varPhi,\tilde{\nabla}_i\varPhi \rangle+[(2\Lambda+\kappa'V)u^{-2}]_i\\&
\ \ +\kappa'(\tilde{\triangle}\phi\phi_i+\frac{1}{2}|d\phi|^2_{\tilde{g}, i}),
\end{split}
\end{equation}
since $tr_{\tilde{g}}(\varPhi^{\ast}h)=2\tilde{R}-2\kappa'|d\phi|^2_{\tilde{g}}-(12\Lambda +6\kappa'V)u^{-2}$,  we get $$\frac{1}{2}\langle\tilde{\triangle}\varPhi,\tilde{\nabla}_i\varPhi \rangle-[(\Lambda+\frac{\kappa'}{2}V)u^{-2}]_i +\frac{1}{2}u^{-2}\kappa' (V\circ\phi)_i=0.$$ Together  with  ii) in Lemma \ref{l2.2}, we have 
$$
\langle\omega,B\rangle E+\langle\omega, E\rangle B=0.
$$
This gives $2\lambda \langle\omega, E\rangle E=0,$ which implies  $\lambda=0$. Hence $B\equiv 0$ near $P$, from the last equation in (\ref{2.19}), we know $\langle \omega, E\rangle\equiv 0$ near $P$, which is a contradiction.   If $ \langle\omega, B\rangle\neq 0$, a contradiction can be derived  by the same argument. 
\end{pf}

Now we can apply the standard Bochner formula(for a map not necessarily harmonic, see \cite{EL}):  
\begin{equation}\label{2.31}
\hat{\triangle} e(\varPhi)= 2\langle \hat{\nabla} \varPhi, \hat{\nabla}\hat{\triangle}\varPhi \rangle_{\hat{g}}+2|\hat{\nabla}_{\alpha\beta}\varPhi|^2+2\langle\hat{R}ic, \varPhi^{\ast}h \rangle_{\hat{g}} -2R_{abcd}\varPhi^{a}_{\alpha}\varPhi^{b}_{\beta}\varPhi^{c}_{\gamma}\varPhi^{d}_{\delta}\hat{g}^{\beta\delta} \hat{g}^{\alpha\gamma}
\end{equation}
where \begin{equation} \label{2.32}
e(\varPhi)= {\hat{g}}^{\alpha\beta}(\varPhi^{\ast}g_{H})_{\alpha\beta}=4|\nabla\log u|^2+u^{-4}|\omega|^2-2\kappa u^{-2}|E|^2-2\kappa u^{-2}|B|^2.
\end{equation}  

By using (\ref{2.25}), Lemma \ref{2.3} and the following formula, 
$$
\hat{\nabla}_{\alpha}\hat{\triangle}\varPhi^a=\frac{\partial}{\partial x^{\alpha}}\hat{\triangle}\varPhi^a+\Gamma^{a}_{bc}\hat{\triangle}\varPhi^b \varPhi^{c}_{\alpha}, 
$$ one can get 
\begin{equation}\label{2.33}
\begin{split}
\hat{\nabla}_{\alpha}\hat{\triangle}\varPhi^1&=[-(2\Lambda+\kappa'V) u^2]_{\alpha}+\Gamma^{1}_{11}(-2\Lambda-\kappa' V) u^2  \varPhi^{1}_{\alpha}=-\kappa'u^2 (V\circ\phi)_{\alpha} \\
\hat{\nabla}_{\alpha}\hat{\triangle}\varPhi^2&= \Gamma^{2}_{12}(-2\Lambda-\kappa'V) u^2) \varPhi^{2}_{\alpha}=(2\Lambda+\kappa'V) \omega_{\alpha} \\
\hat{\nabla}_{\alpha}\hat{\triangle}\varPhi^3&=\Gamma^{3}_{13}(-2\Lambda-\kappa'V) u^2 E_{\alpha}=(\Lambda+\frac{\kappa'}{2}V) E_{\alpha}\\
\hat{\nabla}_{\alpha}\hat{\triangle}\varPhi^4&=(\Lambda+\frac{\kappa'}{2}V) B_{\alpha}.
\end{split}
\end{equation}
Hence \begin{equation}\label{2.34}
\begin{split}\langle \hat{\nabla} \varPhi, \hat{\nabla}\hat{\triangle}\varPhi \rangle_{\hat{g}}=& (2\Lambda+\kappa'V)[|\omega|^2u^{-4}-\kappa u^{-2}(|E|^2+|B|^2)]\\&-2\kappa'\langle \nabla (V\circ \phi), \nabla \log u \rangle.
\end{split}
\end{equation}

 Now we compute the term $|\hat{\nabla}_{\alpha\beta}\varPhi|^2$ at a fixed point $P=(x',t')$. By definition,  $$|\hat{\nabla}_{\alpha\beta}\varPhi|^2=u^{-4}|(\hat{\nabla}_{\alpha\beta}\varPhi)^1|^2_{\hat{g}}+u^{-4}|(\hat{\nabla}_{\alpha\beta}\varPhi)^2|^2_{\hat{g}}-2\kappa u^{-2} |(\hat{\nabla}_{\alpha\beta}\varPhi)^3|^2_{\hat{g}}-2\kappa u^{-2} |(\hat{\nabla}_{\alpha\beta}\varPhi)^4|^2_{\hat{g}}.$$

 Let $\{x^{i}\}$ be a spatial normal coordinate system  around  the fixed point ${x}'$.     
Let $F_{0}=u^{-1}\frac{\partial}{\partial t}$, $F_{i}=\frac{\partial}{\partial x^i}-\theta_i\frac{\partial}{\partial t}$,  then $\{F_{\alpha}\}$ is an orthonormal basis of $\hat{g}$ at $({x}',t')$, hence 
\begin{equation}\label{2.35}
\begin{split}
|(\hat{\nabla}_{\alpha\beta}\varPhi)^a|^2_{\hat{g}}=\sum_{\alpha,\beta}[(\hat{\nabla}^2\varPhi)^a(F_{\alpha},F_{\beta})]^2.
\end{split}
\end{equation}

Combining   (\ref{2.28}) and (\ref{2.5}),  we get 
\begin{equation}\label{2.36}
\begin{split}
u^{-4}|(\hat{\nabla}_{\alpha\beta}\Phi)^1|^2_{\hat{g}}=& 4|\nabla \log u|^4+2|u^{-2}\omega\wedge d\log u|^2\\& +|2\nabla_{ij}\log u+u^{-4}\omega_i\omega_j-\kappa u^{-2}(E_iE_j+B_iB_j)|^2\\
u^{-4}|(\hat{\nabla}_{\alpha\beta}\Phi)^2|^2_{\hat{g}}
=& \langle d \log u,u^{-2}\omega\rangle^2
+u^{-4}|\nabla_i\omega_j-2\omega_i(\log u)_j-2\omega_j(\log u)_i|^2.
\end{split}
\end{equation}
Similarly, one can compute 
\begin{equation}\label{2.37}
\begin{split}
u^{-2}|(\hat{\nabla}_{\alpha\beta}\Phi)^3|^2_{\hat{g}}=&\langle d\log u,u^{-1}E\rangle^2+2|\frac{1}{2}u^{-2}\omega\wedge u^{-1}E|^2\\& +u^{-2}|\nabla_iE_j-E_i(\log u)_j-E_j(\log u)_i|^2\\
u^{-2}|(\hat{\nabla}_{\alpha\beta}\Phi)^4|^2_{\hat{g}}=&\langle d\log u,u^{-1}B\rangle^2+2|\frac{1}{2}u^{-2}\omega\wedge u^{-1}B|^2\\& +u^{-2}|\nabla_iB_j-B_i(\log u)_j-B_j(\log u)_i|^2.
\end{split}
\end{equation}

To compute the fourth   term on the right hand side of (\ref{2.31}),   we need 
\begin{equation}\begin{split}\label{2.38}
   &R_{abcd}\varPhi^{a}_{\alpha}\varPhi^{b}_{\beta}\varPhi^{c}_{\gamma}\varPhi^{d}_{\delta}\hat{g}^{\alpha\gamma} \hat{g}^{\beta\delta}=R_{abcd} \langle \nabla y_a, \nabla y_c\rangle \langle \nabla y_b, \nabla y_d\rangle\\
   &=2R_{1212}(|\nabla u^2|^2|\omega|^2-\langle \nabla u^2,\omega\rangle^2)+2R_{1313} (|\nabla u^2|^2|E|^2-\langle \nabla u^2,E\rangle^2)\\
   & \ \ +2R_{1414}(|\nabla u^2|^2|B|^2-\langle \nabla u^2,B\rangle^2)+2R_{2323}(|\omega|^2|E|^2-\langle\omega,E\rangle^2)\\&\ \  +2R_{2424} (|\omega|^2|B|^2-\langle\omega,B\rangle^2)
   +2R_{3434} (|E|^2|B|^2-\langle E,B\rangle^2),
   \end{split}
\end{equation}
where we have used Lemma \ref{l2.1}. 

From Lemma \ref{l2.1},  we know 
\begin{equation}\begin{split}\label{2.39}
& R_{1212}=-u^{-8}, R_{1313}=R_{1414}=\frac{\kappa}{2}u^{-6}\\ 
& R_{2323}=R_{2424}=\kappa u^{-6}, R_{3434}=-\kappa^2 u^{-4}.
   \end{split}
\end{equation}
So
\begin{equation}\begin{split}\label{2.40}
   &-R_{abcd}\varPhi^{a}_{\alpha}\varPhi^{b}_{\beta}\varPhi^{c}_{\gamma}\varPhi^{d}_{\delta}\hat{g}^{\alpha\gamma} \hat{g}^{\beta\delta}\\ &= 8|d\log u\wedge u^{-2}\omega|^2-4\kappa |d\log u\wedge u^{-1}E|^2-4\kappa |d\log u\wedge u^{-1}B|^2\\& \ \ -2\kappa |u^{-2}\omega\wedge u^{-1}E|^2-2\kappa |u^{-2}\omega\wedge u^{-1}B|^2+2\kappa^2 |u^{-1}E\wedge u^{-1}B|^2.
   \end{split}
\end{equation}
To compute the  third  term on the right hand side of (\ref{2.31}), we rewrite (see (\ref{2.11}))
\begin{equation}\nonumber
\begin{split}
\bar{R}ic(e_k,e_l)=&-\frac{\kappa}{2}u^{-2}(|E|^2+|B|^2)g_{kl}+\kappa u^{-2}(|E|^2 g_{kl}-E_kE_l)\\ &+\kappa u^{-2}
(|B|^2 g_{kl}-B_kB_l)+(\Lambda+\frac{\kappa'}{2}V) g_{kl}+\kappa' \phi_{k}\phi_l,
\end{split}
\end{equation}
thus 
\begin{equation}
\begin{split}\label{2.41}
&\langle\hat{R}ic, \varPhi^{\ast} g_{H}\rangle_{\hat{g}}\\
&= [u^{-4}(\omega^{i}\omega^{j}-|\omega|^2g^{ij})+\bar{R}ic(e_k,e_l)g^{ik} g^{jl}] \\ & \ \ \ \times [u^{-4}\omega_{i}\omega_{j}+4u^{-2}{u_iu_j}-2\kappa u^{-2}E_iE_j-2\kappa u^{-2} B_iB_j]\\
&=[\Lambda+\frac{\kappa'}{2}V-\frac{\kappa}{2}u^{-2}(|E|^2+|B|^2)]e(\varPhi)-4|u^{-2}\omega\wedge d\log u|^2\\& \ \ \ +3\kappa |u^{-2}\omega \wedge u^{-1}E|^2+3\kappa |u^{-2}\omega\wedge u^{-1}B|^2+4\kappa|u^{-1}E\wedge d\log u|^2\\& \ \ \ +4\kappa|u^{-1}B\wedge d\log u|^2-4\kappa^2|u^{-1}E\wedge u^{-1}B|^2+\kappa'\langle d\phi\otimes d\phi, \varPhi^{\ast} g_{H}\rangle.
\end{split}
\end{equation}

Note  that $E\wedge B=0$ holds by condition (\ref{2.20}). Combining 
\begin{equation}\label{2.42}
\hat{\triangle} (\frac{\kappa'}{2}V\circ \phi)=\frac{\kappa'}{2} Hess(V)_{ab}\phi^a_{\alpha}\phi^b_{\beta}\hat{g}^{\alpha\beta}+\frac{\kappa'}{4}|\nabla V|^2
\end{equation}

   and  (\ref{2.31}) (\ref{2.34}) (\ref{2.36}) (\ref{2.37})(\ref{2.40})(\ref{2.41}), we have

\begin{equation}\label{2.43}
\begin{split}
&\hat{\triangle}(\frac{1}{2}e(\varPhi)+\frac{\kappa'}{2}V)\\ &=  4|\nabla \log u|^4 +|2\nabla_{ij}\log u +u^{-4}\omega_i\omega_j-\kappa u^{-2}(E_iE_j+B_iB_j)|^2\\& \ \ +u^{-4}|\nabla_i\omega_j-2\omega_i(\log u)_j-2\omega_j(\log u)_i|^2\\
&\ \ -2\kappa u^{-2}|\nabla_iE_j-E_i(\log u)_j-E_j(\log u)_i|^2\\
&\ \ -2\kappa u^{-2}|\nabla_iB_j-B_i(\log u)_j-B_j(\log u)_i|^2\\& \ \ +\langle d \log u,u^{-2}\omega\rangle^2-2\kappa \langle d\log u,u^{-1}E\rangle^2-2\kappa \langle d\log u,u^{-1}B\rangle^2\\
&\ \ +6|u^{-2}\omega\wedge d\log u|^2+[\Lambda+\frac{\kappa'}{2}V-\frac{\kappa}{2}u^{-2}(|E|^2+|B|^2)]e(\varPhi)\\&\ \ +(2\Lambda +\kappa' V)[|\omega|^2u^{-4}-\kappa u^{-2}(|E|^2+|B|^2]+I
 \end{split} 
 \end{equation}
 where \begin{equation} \label{2.44}
\begin{split}
I&= \frac{\kappa'}{2}[Hess(V)_{ab} \phi^a_\alpha \phi^a_\beta\hat{g}^{\alpha\beta}+\frac{1}{2}|\nabla V|^2]+\kappa'\langle d\phi\otimes d\phi, \varPhi^{\ast} g_{H}\rangle\\&\ \ \ \ -2\kappa'\langle \nabla (V\circ \phi), \nabla \log u \rangle\\
& \geq \kappa'|\frac{1}{2}\nabla V-2d\phi(\nabla\log u)|^2\\ &  \geq 0
\end{split}
\end{equation}
for $\kappa'\geq 0$.

\begin{lem} \label{l2.4} For any stationary solution to (\ref{EMKG}) and (\ref{1.4}) with $\Lambda\geq 0, \kappa\leq 0, \kappa'\geq 0$,  under the assumption (\ref{2.20}),  we have 
\begin{equation}\label{2.45}
\begin{split}
& \hat{\triangle}(\frac{1}{2}e(\varPhi)+\frac{\kappa'}{2}V)= \triangle (\frac{1}{2}e(\varPhi)+\frac{\kappa'}{2}V)
+\langle\nabla \log u, \nabla(\frac{1}{2}e(\varPhi)+\frac{\kappa'}{2}V)\rangle\\& \geq   4|\nabla \log u|^4 +|2\nabla_{ij}\log u +u^{-4}\omega_i\omega_j-\kappa u^{-2}(E_iE_j+B_iB_j)|^2\\&\ \ \ +[\Lambda+\frac{\kappa'}{2}V-\frac{\kappa}{2}u^{-2}(|E|^2+|B|^2)]e(\varPhi),\end{split} 
 \end{equation}
 where \begin{equation}\label{2.46}
 \frac{1}{2}e(\varPhi)+\frac{\kappa'}{2}V=2|\nabla\log u|^2+\frac{1}{2}u^{-4}|\omega|^2-\kappa u^{-2}|E|^2-\kappa u^{-2}|B|^2+\frac{\kappa'}{2}V.
 \end{equation}
 \end{lem}
The most notable feature   in formula  (\ref{2.45}) is that each term in the right hand side of  (\ref{2.45}) is nonnegative if $\Lambda\geq 0,\kappa\leq0$ and $\kappa'\geq 0$.

\subsection{Proof of Theorems \ref{N1.1} and \ref{N1.2}.}
\begin{pf} of Theorem \ref{N1.1}.  The idea of the proof is analogous to  Theorem 5.3 in  \cite{Chen1}. We assume $\partial \hat{B}(x_0,a)\neq \phi$.

Let $h(x)=2|\nabla \log u|^2(x)+\frac{1}{2}u^{-4}|\omega|^2(x)-\kappa u^{-2}|E|^2-\kappa u^{-2}|B|^2+\frac{\kappa'}{2}V+\Lambda$ be the quantity in (\ref{2.46})(up to a constant $\Lambda$), $f(x)=h(x)d_{\hat{g}}^2(x,\partial \hat{B}(x_0,a))$, and 
 $\bar{x}\in \hat{B}(x_0,a)$ such that $f(\bar{x})=\sup_{x\in \hat{B}(x_0,a)}f(x)$.
 
 For any fixed $0<\delta<1$,  we have two cases, 1) $h(\bar{x})\leq  {\delta} \frac{m^2}{\hbar^2}$, 2)  $h(\bar{x})> {\delta} \frac{m^2}{\hbar^2}$. 

For case 1), for any $x\in \hat{B}(x_0,\frac{3a}{4})$, we have 
$(\frac{a}{4})^2 h(x) \leq f(\bar{x}) \leq \delta \frac{m^2}{\hbar^2} a^2$, which implies $$\sup_{x\in \hat{B}(x_0,\frac{3a}{4})}h(x)\leq 16\delta  \frac{m^2}{\hbar^2}. $$ 

 So we may assume case 2) always holds, and we will show that $f(\bar{x})\leq C_{\delta}$ for some constant $C_\delta$ depending only on $\delta$. 
 Note that 
 \begin{equation}\label{2.47}
 \kappa' |\nabla (V\circ \phi)|\leq 2\kappa' \frac{m^2}{\hbar^2}|d\phi||\phi| \leq \kappa' \frac{m^2}{\hbar^2}V +\kappa' |d\phi|^2.  \end{equation}

 We will argue by contradiction.  Suppose there are  a sequence of 4-Lorentzian manifolds $(M_l, \bar{g}_{l})$ satisfying the  equations (\ref{EMKG}) (\ref{1.4}),  and a sequence of $\hat{g}_l$-balls $\hat{B}(x_l,a_l)\subset  M_{l}$  with compact closure such that $f(\bar{x}_l)\rightarrow \infty$  as $l\rightarrow \infty$,  where 
\begin{equation}\label{2.48}
\begin{split}
& f(\bar{x}_l)=\sup_{x\in \hat{B}(x_l,a_l)} h_{l}(x)d^2_{\hat{g}_l}(x,\partial \hat{B}(x_l,a_l))\\& h_{l}(x)=2|\nabla \log u_l|^2+\frac{1}{2}u_l^{-4}|\omega_l|^2-\kappa u_l^{-2}|E_l|^2-\kappa_l u_l^{-2}|B_l|^2+\Lambda_{l}+\frac{\kappa'_l}{2}V_l.
\end{split}
\end{equation}
Scaling $u_l$ and $\bar{g}_l$ by $u_l(\bar{x}_l)^{-1}$ and $h_l(\bar{x}_l)$ respectively, one can assume  $u_l(\bar{x}_l)=1$ and  $h_l(\bar{x}_l)=1$.  We still use the same notations $u_l,$ $\omega_l$, $E_l$, $B_l$, $h_l$, $\bar{g}_l$, etc.  to denote the corresponding scaled quantities.  
For any fixed $ 0<\epsilon<1$,  any $x\in \hat{B}(x_0,a_l)$ with $$d_{\hat{g}_l}(x,\bar{x}_l)\leq \epsilon f(\bar{x}_l)^{\frac{1}{2}} h_l^{-\frac{1}{2}}(\bar{x}_l)=\epsilon d_{\hat{g}_l}(\bar{x}_l,\partial \hat{B}(x_0,a_l)),$$
 we have $d_{\hat{g}_l}(x, \partial \hat{B}(x_0,a_l))\geq (1-\epsilon) d_{\hat{g}_l}(\bar{x}, \partial \hat{B}(x_0,a_l))$, hence 
 \begin{equation} \label{2.49}
 h_l(x)\leq \frac{1}{(1-\epsilon)^2} h_l(\bar{x}_l).
  \end{equation}
  
 If there is no ambiguity, we can also omit the subscript $l$ from $h_l$, $u_l$, etc.   It can be shown from (\ref{2.49}): 
  \begin{equation} \label{2.50}
 h(x)\leq 4\ \ \  \text{on}  \  \hat{B}(\bar{x}, D). 
  \end{equation}
  where 
       $D=\frac{1}{2}\sqrt{f(\bar{x})}.$
         In the following, we  estimate $\kappa' |d\phi|^2$ on $\hat{B}(\bar{x},\frac{D}{2})$.   
From (\ref{2.10})(\ref{2.11}) (\ref{2.14}) and (\ref{2.4}), we know \begin{equation} \label{2.51}
\begin{split}
& u^{-2}\hat{R}ic(X,X)\geq \frac{\kappa}{2}u^{-2}(|E|^2+|B|^2)\\
& u^{-1}\hat{R}ic(X,e_j)=\kappa u^{-2} i_{e_j}(\ast(E\wedge B))
\\& \hat{R}ic({e_i,e_j}) \geq \frac{\kappa}{2}u^{-2}(|E|^2+|B|^2)g_{ij} -u^{-4}|\omega|^2g_{ij}.
\end{split}
\end{equation}
  Together with  (\ref{2.50}) and  (\ref{2.48}),  we get
  \begin{equation}\label{2.52} \hat{R}ic\geq -2h\hat{g}\geq -8\hat{g} \ \  \text{on}  \  \hat{B}(\bar{x}, D).
  \end{equation}
  Note that we have the Bochner formula: 
\begin{equation}\label{2.53}
\begin{split}
\hat{\triangle}\kappa' |d\phi|^2&=2\kappa' \langle\hat{R}ic,d\phi\otimes d\phi\rangle
 +2\kappa' |\hat{\nabla}d\phi|^2+2\kappa' \frac{m^2}{\hbar^2}|d\phi|^2\\&\geq -2\kappa' |d\phi\wedge u^{-2}\omega|^2+\kappa \kappa' u^{-2}(|B|^2+|E|^2)|d\phi|^2+\frac{2}{3}\kappa'^2|d\phi|^4.
\end{split}
\end{equation}

  Let $\xi:[0,\infty) \rightarrow [0,1]$ be a fixed nonnegative smooth non-increasing function such that $\xi=1$ on $[0,\frac{1}{2}]$ and $\xi=0$ on $ [1,\infty)$. Consider the  function $L(x)=\xi(\frac{d_{\hat{g}}(\bar{x},x)}{D})(\kappa' |d\phi|^2)$ which is nonnegative and vanishes on $\partial \hat{B}(\bar{x},D)$. So $L(x)$ assumes its maximum at some point $\hat{x}\in B(\bar{x}, D)$.  We temporarily assume $\hat{x}$ is a smooth point of $d_{\hat{g}}(\bar{x},\cdot)$.  If $d_{\hat{g}}(\hat{x},\bar{x})\geq \frac{D}{2}$,  by (\ref{2.53}) and  Laplacian comparison theorem (see Corollary 1.2 in \cite{SY} ), one can show 
  \begin{equation}\label{2.54}
  \hat{\triangle} d_{\hat{g}}(\bar{x},\cdot)\mid_{x=\hat{x}}\leq  \frac{3}{d_{\hat{g}}(\bar{x},\hat{g})}+3 \sqrt{\frac{8}{3}}\leq \frac{6}{D}+6,
  \end{equation}
  hence  
  \begin{equation}\label{2.55}
  \begin{split}
  0&\geq \hat{\triangle}L(x)\mid_{x=\hat{x}}\\&=\xi \hat{\triangle} (\kappa'|d\phi|^2)-\frac{2}{D^2}\frac{(\xi^{'})^2}{\xi}(\kappa'|d\phi|^2)+(\kappa'|d\phi|^2)(\frac{1}{D^2}\xi^{''}+\frac{1}{D}\xi^{'}\hat{\triangle}{d_{\hat{g}}})\\
  & \geq \frac{2}{3}(\kappa'|d\phi|^2)^2\xi-[|\kappa|u^{-2}(|B|^2+|E|^2)+2u^{-4}|\omega|^2+CD^{-1}]\kappa'|d\phi|^2\\
  &\geq \frac{2}{3}(\kappa'|d\phi|^2)^2\xi-16 \kappa'|d\phi|^2-{C}D^{-1}\kappa'|d\phi|^2,
  \end{split}
  \end{equation}
  where we have used (\ref{2.48}) and (\ref{2.50}). 
   If $\hat{x}$ lies in the cut locus of $\bar{x}$, we can 
   use a standard support function technique as in \cite{Y} (or Theorem 3.1 in \cite{SY}) to prove that (\ref{2.55}) still holds. 

  Multiplying  both sides of (\ref{2.55}) by $\xi(\hat{x})$,  we find 
  $$
  L(\hat{x})^2-CL(\hat{x})\leq 0.
  $$
  which implies $L(\hat{x})\leq C$. Therefore,  
  \begin{equation}\label{2.56}
  \sup_{x\in \hat{B}(\bar{x},\frac{D}{2})}\kappa'|d\phi|^2\leq C
  \end{equation}
   where $C$ is a universal constant.   
  Combining (\ref{2.50}) and (\ref{2.56}),  the estimate (\ref{2.47}) becomes 

  \begin{equation}\label{2.57}
 \kappa' |\nabla (V\circ \phi)|\leq  8 \delta^{-1} +C, \ \ \ \ \text{on} \ \  \hat{B}(\bar{x},\frac{D}{2}),
\end{equation}
because the quantity  $\frac{m^2}{\hbar^2}$ in (\ref{2.47}) has been changed to $\frac{m^2}{\hbar^2}h(\bar{x})^{-1}\leq \delta^{-1}$ after scaling. 

From (\ref{2.16}) and (\ref{2.50}) (\ref{2.57}),    the Ricci curvature of $\tilde{g}$ is uniformally bounded (independent of $l$)  on $\hat{B}(\bar{x},1)$. As in the proof of Theorem 5.3 in \cite{Chen1},  one can use  the horizontal exponential map (w.r.t. metric $u_l^2 \bar{g}_l$) to pull back $\tilde{g}_l$, $\omega_l$ , $E_l$, $B_l$, $u_l$ to horizontal tangent space,  moreover, using \cite{JK}, one can construct a $\tilde{g}$-harmonic coordinate system  $\{z^{i}\}$  (of uniform size $\{|z|<\delta_1\}$, $\delta_1$ independent of $l$) on the horizontal tangent space around $\bar{x}_l$, so that the $C^{1,\alpha}$-norm (w.r.t. $\{z^i\}$ coordinates) of $\tilde{g}$ is uniformally bounded, see Theorem 5.3  (5.31)  in  \cite{Chen1}.

Note that   $u_l, \omega_l, E_l, B_l, \phi_l$ satisfy the following elliptic type equations (see (\ref{2.6}) (\ref{2.19}) (\ref{2.21})): 
\begin{equation} \label{2.58}
\begin{split}
 u^2\tilde{\triangle}\log u&=-\frac{1}{2}u^{-4}|\omega|^2+\kappa (\frac{1}{2} u^{-2}|E|^2+\frac{1}{2}u^{-2}|B|^2)-(\Lambda+\frac{\kappa'}{2}V)\\
    u^2\tilde{\triangle}\phi_3&=2\langle d\log u, E\rangle+u^{-2}\langle \omega,B\rangle\\
     u^2\tilde{\triangle}\phi_4&=2\langle d\log u, B\rangle+u^{-2}\langle \omega,E\rangle\\
     u^{2}\tilde{g}^{kl}\tilde{\nabla}_k\omega_l&=4\tilde{g}^{kl}\omega_k\nabla_l \log u\\
     d\omega&=\pm 2\kappa \ast (E\wedge B)\\
     u^2\tilde{\triangle}\phi&=\frac{m^2}{\hbar^2}\phi
\end{split}
\end{equation}
where $0\leq \frac{m^2}{\hbar^2}\leq \delta^{-1}$. In the following, we say the $C^{k,\alpha}$-norms  of  quantities  $F_{l}$ are  uniformally bounded, if for any $0<\delta_2< \delta_1$, there is a constant $C_{\delta_2,k, \alpha}$( independent of $l$) such that $|F_l|_{C^{k,\alpha}(\{|z|\leq\delta_2\})}\leq C_{\delta_2,k, \alpha}$.  From the boundedness of $h_l$ , using   equations  1-5 in (\ref{2.58}) and elliptic regularity,  one can show that the $C^{1,\alpha}$-norms   of $\log u_l$ are  uniformally bounded, hence the $C^{1,\alpha}$-norms of  $\sqrt{|\kappa_l|}E_l,\sqrt{|\kappa_l|}B_l,\omega_l$ are also uniformally bounded.  Combining (\ref{2.57}) and the first equation in (\ref{2.58}),  we know the  $C^{2,\alpha}$-norms  of $\log u_l$  are  uniformally bounded.  Applying equations 2-5 in (\ref{2.58}) again,  the  $C^{2,\alpha}$-norms of $\sqrt{|\kappa_l|}E_l,\sqrt{|\kappa_l|}B_l,\omega_l$ are also uniformally bounded.

By differentiating the last equation of (\ref{2.58}), we get 
\begin{equation} \label{2.59}
\begin{split}
      \tilde{\triangle}(\sqrt{\kappa'}\phi_i)&=\frac{m^2}{\hbar^2}\sqrt{\kappa'}(u^{-2}\phi)_i+\tilde{R}_{ik}\sqrt{\kappa'}\phi_j\tilde{g}^{kj}.
\end{split}
\end{equation}
Since $\sqrt{\kappa'}\phi_i$, $\tilde{R}_{ij}$  are  uniformally bounded and  $\frac{m^2}{\hbar^2}\leq \delta^{-1}$, $\kappa' \frac{m^2}{\hbar^2}|\phi|^2\leq 8$, by $L^p$-estimate for elliptic equations,  one can show  the $C^{1,\alpha}$-norms  of $\sqrt{\kappa'}\phi_{i}$ are  uniformally bounded. So the $C^{3,\alpha}$-norms  of $\tilde{g}_l$ are  uniformally bounded (using the   harmonic coordinates $\{z^{\i}\}$ in (\ref{2.16})). By repeating the above arguments, one can show that for any $k\in Z_{+}$, the $C^{k,\alpha}$-norms of $\log u_l$, $\omega_l$, $\sqrt{-\kappa_l} E_l$, $\sqrt{-\kappa_l}  B_l$, $\sqrt{\kappa'} d\phi_l$,  $\kappa'_l V_l$ are uniformally bounded. 
So we can extract a smooth convergent subsequence so that the limit $u^{\infty}, \tilde{g}^{\infty},  (\sqrt{|\kappa|}E)^{\infty},(\sqrt{|\kappa|}B)^{\infty},  (\sqrt{\kappa'} d\phi)^{\infty}, (\kappa'V)^{\infty}\in C^{\infty}(|z|<\delta_1)$.   Note that the condition (\ref{1.7}) imply that the limit must satisfy  $ (\sqrt{|\kappa|}E)^{\infty}\wedge (\sqrt{|\kappa|}B)^{\infty}=0$. Hence equation (\ref{2.45}) holds for the limit.  The smooth  convergence of $h_{l}$ and (\ref{2.49})  imply that $h^{\infty}(x)$ will achieve its maximum ($=1$) at the origin $z^i=0$.    Since each term in the right hand side of (\ref{2.45}) is nonnegative,  we can apply the strong maximum principle for the limit. 
This implies that  $h^{\infty}\equiv const.=1$  and each term on the right hand side of (\ref{2.45}) vanishes. From $|\nabla\log u^{\infty}|^4=0$, we know $u^{\infty}\equiv const.=1$. Combining  with   $|[2{\nabla}_{ij}\log u+u^{-4}\omega_i\omega_j-\kappa u^{-2}(E_iE_j+B_iB_j)]^{\infty}|^2=0$, we find $(\sqrt{\kappa}E)^{\infty}=(\sqrt{\kappa}B)^{\infty}=\omega^{\infty}=0$. From the expression of $h$, we find  $[\Lambda+\frac{1}{2} \kappa' V]^{\infty}\equiv 1$. The first equation in (\ref{2.58}) for the limit will give a contradiction:  $0=(u^{\infty})^2\tilde{\triangle}\log u^{\infty}=-(\Lambda+\frac{\kappa'}{2}V)^{\infty}=-1$.  
In conclusion, we have proved
\begin{equation}\sup_{\hat{B}(x_0,\frac{3a}{4})} h(x)\leq  16 \delta \frac{m^2}{\hbar^2}+C_{\delta} a^{-2}.
\end{equation} 
From the above proof (on the regularity of $\log u$),    $\sup_{\hat{B}(x_0,\frac{5a}{8})}|u^{-1}\nabla^2 u|+u^{-2}|\nabla \omega|$ can also be bounded  by $ C \delta \frac{m^2}{\hbar^2}+CC_{\delta} a^{-2}$.  Now combining   (2.10) in \cite{Chen1} and formula  (\ref{2.2}), we know $\sup_{\hat{B}(x_0,\frac{a}{2})}|Rm|_{\hat{g}}\leq C \delta \frac{m^2}{\hbar^2}+CC_{\delta} a^{-2}$.  The estimate (\ref{1.8}) follows by redefining  the constants.   

If  $\partial \hat{B}(x_0,a)$ is empty,  by using  strong maximum principle on  equation (\ref{2.45}), we conclude that $Rm\equiv 0$, which fulfills  the estimate (\ref{1.8}). 

\end{pf}
\begin{pf} of Theorem \ref{N1.2}. By Theorem 3.3 in \cite{Chen1}, $(M, \hat{g})$ is geodesically complete. First, let $a\rightarrow\infty$ in (\ref{1.8}),  we find $|Rm|_{\hat{g}}\leq \delta$.  Let $\delta\rightarrow0$,  we find $Rm\equiv 0$, i.e., $(M,g_{M})$ is flat. 
\end{pf}
\section{Higher dimensional static  fields}

 In this section, we generalize  our estimate to general  dimensions.  We may assume the scalar fields can  take their values on some manifold. Let $(M,g_{M})$ be a spacetime of dimension $n+1$, $(W, g_{W})$ a Riemannian manifold equipped with a Riemannian metric $g_{W}$, $V:W\rightarrow R$ be a fixed  function.   We consider the  (uncharged) Einstein-Maxwell-Klein-Gordon   equation, which  is the following system of differential equations on the Lorentzian metric $g_{M}$, a 2-form $F$ on $M$,  and a map (scalar field) $\phi: M\rightarrow W$: 
 \begin{equation}\label{3.1}
 \begin{split}
&Ric(g_{M})-\frac{R}{2}g_M+\Lambda g_M=\kappa(F_{\alpha\gamma}F_{\beta\delta}g^{\gamma\delta}-\frac{1}{4}|F|^2g_{\alpha\beta})\\& \ \ \ \ \ \ \ \ \ \ \  \ \ \ \ \ +\kappa'[\phi^{\ast}g_{W}-\frac{1}{2}(|d\phi|^2+V(\phi))g_{M}]\\
&dF=d\ast F=0\\
&\triangle_{g_{M},g_{W}} \phi =\frac{1}{2} (\nabla V)\circ \phi
\end{split}
\end{equation}
where $\Lambda, \kappa,\kappa'$ are  constants, $\triangle_{g_{M},g_{W}} \phi $ is the harmonic map Laplacian of $\phi$. The third  equation in (\ref{3.1}) may be regarded as a generalized Klein-Gordon equation.   Let  $X$ be  a   timelike Killing field on $(M,g_{M})$, as before,  we say the solution $(g_{M},F,\phi)$ to (\ref{3.1})  is inheriting, if $F$ and $\phi$ satisfies 
\begin{equation}\label{3.2}
\begin{split} \mathcal{L}_X(F)=0, \ \ d\phi(X)=0.
 \end{split}
 \end{equation}

  \begin{thm} \label{N3.1}  Let $(M,g_{M})$ be a static spacetime of dimension $n+1\geq 4$ with a timelike Killing field $X$, $\phi: M\rightarrow W$  a map. We assume \\
  i) $(g_{M},F,\phi)$   satisfies   the   equations (\ref{3.1})(\ref{3.2})
 with  $\Lambda\geq 0, \kappa \leq 0, \kappa'\geq 0$;\\
   ii)  the  magnetic field $B=i_{X}\ast  F$ vanishes;\\
 iii)  $(W,g_{W})$ has nonpositive sectional curvature;\\
 iv) the function $V:W\rightarrow R$ is nonnegative and convex, i.e., $Hess(V)\geq 0$. \\
   Let $\hat{B}(x_0, a)$ be a $\hat{g}-$metric ball centered at $x_0$ of radius $a>0$ with compact  closure  in $M$.  Then there is a universal constant $C>0$  such that 
 \begin{equation} \label{3.3}
\begin{split}
\sup_{x\in \hat{B}(x_0,\frac{a}{2})}|\nabla \log u|^2_{\hat{g}}+|\kappa| u^{-2}|E|^2+\kappa' n^{-1}|d\phi|^2 \leq Cn a^{-2},
\end{split}
\end{equation}
where $E=i_{X} F$ is the electric field. 
\end{thm}
Before the proof of Theorem \ref{N3.1}, we have to mention  that under the condition  that the spacetime is diffeomorphic to $N^n \times R$ so that (\ref{2.1}) holds globally on $N^n \times R$,  analogous result as in Theorem \ref{N3.1} on static Einstein-scalar equation  has been obtained in \cite{Re}. 

  Taking trace on the first equation of (\ref{3.1}), we get  $$
R=2\frac{n+1}{n-1}(\Lambda+\frac{\kappa'}{2}V)+\kappa\frac{n-3}{2(n-1)}|F|^2+\kappa'|d\phi|^2.
$$ Substituting it into (\ref{3.1}),  we find  \begin{equation}\label{3.4}
R_{\alpha\beta}=\kappa F_{\alpha\gamma}F_{\beta\delta}g^{\gamma\delta}+\kappa' (d\phi\otimes d\phi)_{\alpha\beta}+\frac{1}{n-1}(2\Lambda+\kappa'V-\frac{\kappa}{2}|F|^2)g_{\alpha\beta},
\end{equation}
where $|F|^2=F_{\xi\gamma}F_{\eta\delta}g^{\xi\eta}g^{\gamma\delta}=2u^{-2}(|B|^2-|E|^2)=-2u^{-2}|E|^2$.  
 Let  $\{x^{\alpha}\}$ be a local coordinate system so that (\ref{2.1}) holds.  From (\ref{3.4}) and (\ref{2.13}),  one can show 
\begin{equation}\label{3.5}
\begin{split}
\bar{R}ic(X,X)=&\kappa \frac{n-2}{n-1}|E|^2-\frac{2\Lambda+\kappa'V}{n-1}u^2\\
\bar{R}ic(X,e_j)=&0\\
\bar{R}ic(e_i,e_j)=&\kappa (-u^{-2}E_iE_j+\frac{1}{n-1}u^{-2}|E|^2 g_{ij})+\frac{2\Lambda+\kappa'V}{n-1}g_{ij}\\&+\kappa'\langle d\phi(e_i), d\phi(e_j)\rangle.
\end{split}
\end{equation}

 The function $u$ satisfies 
\begin{equation}\label{3.6}
\begin{split}
\hat{\triangle}{\log u}=u^{-2}\bar{R}ic(X,X)=\kappa\frac{n-2}{n-1}|E|^2u^{-2}-\frac{2\Lambda+\kappa'V}{n-1}.
\end{split}
\end{equation}
The Maxwell equation $dF=d\ast F=0$ is now equivalent to  
\begin{equation}\label{3.7}
\begin{split}
E&=d\phi_3\\
\hat{\triangle}{\phi_3}&=2\langle d \log u, E\rangle
\end{split}
\end{equation}
for some locally defined function $\phi_3$.  We consider a map $\Phi=(y_1,y_2)=(u^2,\phi_3)$ from the  coordinate system $\{x^{\alpha}\}$ to a 2-d upper half space  $H=\{(y_1,y_2): y_1>0\}$ equipped with a metric $g_{H}=y_1^{-2}dy_1^2-4\kappa \frac{n-2}{n-1} y_1^{-1}dy_2^2$. By  similar computations as in Lemma 2.2 (2.28),  one can show:  
\begin{equation}\label{3.8}
\begin{split}
& (\hat{\nabla}_{\alpha\beta} \varPhi)^{1}= 2u^2 \hat{\nabla}_{\alpha\beta}\log u-2\kappa \frac{n-2}{n-1} E_{\alpha} E_{\beta}\\
& (\hat{\nabla}_{\alpha\beta} \varPhi)^{2}=\hat{\nabla}_{\alpha\beta} \phi_3-(\log u)_{\alpha}E_{\beta}-(\log u)_{\beta}E_{\alpha}.
\end{split}
\end{equation}
It  implies that $\hat{\triangle}\Phi=-\frac{4\Lambda+2\kappa'V}{n-1}u^2\frac{\partial}{\partial y_1}$, where we have used (\ref{3.6}) and (\ref{3.7}).

By similar computations as in (\ref{2.34})(\ref{2.36})(\ref{2.37})(\ref{2.40}) and (\ref{2.41}), we have 
\begin{equation}
\begin{split}\label{3.9}
\langle\hat{R}ic, \varPhi^{\ast} g_{H}\rangle_{\hat{g}}&= [\bar{R}ic(e_k,e_l)g^{ik} g^{jl}]  \times [4\frac{u_iu_j}{u^2}-4\kappa\frac{n-2}{n-1} u^{-2}E_iE_j]\\
&=[\frac{2\Lambda+\kappa'V}{n-1}-\kappa \frac{n-2}{n-1} u^{-2}|E|^2]e(\varPhi)+4\kappa|u^{-1}E\wedge d\log u|^2\\
&\ \ \ +4\kappa'|d\phi (\nabla \log u) |^2-4\frac{n-2}{n-1}\kappa'\kappa |d\phi(u^{-1}E^{\#})|^2,
\end{split}
\end{equation}
 \begin{equation}\label{3.10}
\begin{split}
|\hat{\nabla}^2\Phi |^2&= 4|\nabla \log u|^4 +4|\nabla_{ij}\log u-\kappa \frac{n-2}{n-1} u^{-2}E_iE_j|^2\\&\ \ \ -4\frac{n-2}{n-1}\kappa[\langle d\log u,u^{-1}E\rangle^2 +u^{-2}|\nabla_iE_j-E_i(\log u)_j-E_j(\log u)_i|^2],
\end{split}
\end{equation}
\begin{equation}
\begin{split}\label{3.11}
   &-R_{abcd}\varPhi^{a}_{\alpha}\varPhi^{b}_{\beta}\varPhi^{c}_{\gamma}\varPhi^{d}_{\delta}\hat{g}^{\alpha\gamma} \hat{g}^{\beta\delta}= -8\kappa \frac{n-2}{n-1} |d\log u\wedge u^{-1}E|^2,
\end{split}
\end{equation}
\begin{equation}\label{3.12}\langle \hat{\nabla} \varPhi, \hat{\nabla}\hat{\triangle}\varPhi \rangle_{\hat{g}}=-8\frac{n-2}{(n-1)^2}\kappa(\Lambda+\frac{\kappa'}{2}V) u^{-2}|E|^2-\frac{4\kappa'}{n-1}\langle \nabla (V\circ \phi), \nabla \log u\rangle,
\end{equation}
where $E^{\#}$ is the vector field obtained by lifting the indices of $E$. 
Combining (\ref{3.9}) (\ref{3.10}) (\ref{3.11}) (\ref{3.12}) and the  Bochner formula (2.31), we have  
\begin{equation}
\begin{split}\label{3.13}
 &  \hat{\triangle}(|\nabla \log u|^2-\kappa \frac{n-2}{n-1}u^{-2}|E|^2)\\& =2|\nabla \log u|^4+2|\nabla_{ij}\log u-\kappa \frac{n-2}{n-1} u^{-2}E_iE_j|^2-2\frac{n-2}{n-1}\kappa [\langle d\log u,u^{-1}E\rangle^2\\&\ \  \ +|u^{-1}\nabla_iE_j-u^{-1}E_i(\log u)_j-u^{-1}E_j(\log u)_i|^2]-2\frac{n-2}{n-1}\kappa'\kappa  |d\phi(u^{-1}E^{\#})|^2\\
  &\ \ \ -2\kappa \frac{n-3}{n-1}|d\log u\wedge u^{-1}E|^2-\kappa(4\Lambda+2\kappa'V)\frac{n-2}{(n-1)^2}u^{-2}|E|^2\\
  &\ \ \ +(\frac{4\Lambda+2\kappa'V}{n-1}-2\kappa \frac{n-2}{n-1} u^{-2}|E|^2)(|\nabla \log u|^2-\kappa \frac{n-2}{n-1}u^{-2}|E|^2)+I
\end{split}
\end{equation}
where \begin{equation}\begin{split}\nonumber 
I=2{\kappa'}|d\phi(\nabla\log u)|^2-\frac{2\kappa'}{n-1}\langle \nabla (V\circ \phi), \nabla \log u\rangle\geq -\frac{\kappa'}{2(n-1)^2}|\nabla V|^2.
\end{split}
 \end{equation}
 Each term in the right hand side of (\ref{3.13}) is nonnegative except $I$. Note  that the sectional curvature of $(W,g_{W})$ is assumed to be nonpositive and $V$ is convex, we have: \begin{equation}\label{3.14}
\begin{split}
\hat{\triangle}\kappa' |d\phi|^2&=2\kappa' \langle\hat{R}ic,\phi^{\ast}g_{W}\rangle
-2\kappa' R_{abcd}\phi^{a}_{\alpha}\phi^{b}_{\beta}\phi^{c}_{\gamma}\phi^{d}_{\delta}\hat{g}^{\alpha\gamma}\hat{g}^{\beta\delta}\\& \ \ +2\kappa' |\hat{\nabla}d\phi|^2+\kappa' \langle\phi^{\ast}Hess(V), \hat{g}\rangle\\&\geq 2\kappa' \langle\hat{R}ic,\phi^{\ast}g_{W}\rangle+2\kappa' |\hat{\nabla}d\phi|^2.
\end{split}
\end{equation}
Since 
\begin{equation}\label{3.15}
\begin{split}
\langle\hat{R}ic,\phi^{\ast}g_{W}\rangle&=\kappa'\langle d\phi(e_i),d\phi(e_j)\rangle^2+\frac{(2\Lambda+{\kappa'}{V})}{n-1}|d\phi|^2\\& \ \ \ +\kappa[-u^{-2} |d\phi(E^{\#})|^2+\frac{1}{n-1}u^{-2}|E|^2|d\phi|^2]\\
|\hat{\nabla}d\phi|^2& \geq \frac{1}{(n+1)}|\triangle_{\hat{g},g_{W}}\phi|^2=\frac{1}{(n+1)}|\triangle_{\bar{g},g_{W}}\phi|^2=\frac{1}{4(n+1)}|\nabla V|^2,
\end{split}
\end{equation}
we obtain
\begin{equation}\label{3.16}
 \begin{split}
\hat{\triangle}\kappa' |d\phi|^2&\geq \frac{2}{n}\kappa'^2|d\phi|^4+\frac{(4\Lambda+{2\kappa'}{V})}{n-1}\kappa'|d\phi|^2+\frac{\kappa'}{2(n+1)}|\nabla V|^2\\&\ \ +\kappa \kappa' \frac{2}{n-1}u^{-2}|E|^2|d\phi|^2\\
&\geq \frac{1}{8 n}\kappa'^2|d\phi|^4-\frac{8n}{15(n-1)^2}(\kappa u^{-2}|E|^2)^2+\frac{\kappa'}{2(n+1)}|\nabla V|^2.
\end{split}
\end{equation}
Combining (\ref{3.16}) and  (\ref{3.13}),  we get 
\begin{equation}
\begin{split}\label{3.17}
 &  \hat{\triangle}(|\nabla \log u|^2-\kappa \frac{n-2}{n-1}u^{-2}|E|^2+\frac{n+1}{(n-1)^2}\kappa' |d\phi|^2)\\
 & \geq  2|\nabla \log u|^4+\frac{2}{5}\frac{(n-2)^2}{(n-1)^2}\kappa^2(u^{-2}|E|^2)^2+\kappa'^2\frac{n+1}{8n (n-1)^2}|d\phi|^4\\
& \geq \frac{1}{72} (|\nabla \log u|^2-\kappa \frac{n-2}{n-1}u^{-2}|E|^2+\frac{n+1}{(n-1)^2}\kappa' |d\phi|^2)^2,
  \end{split}
\end{equation}
where we have used $2(\frac{n-2}{n-1})^2-\frac{8n}{15}\frac{n+1}{(n-1)^4}\geq \frac{2}{5}(\frac{n-2}{n-1})^2$ and $\frac{1}{8n}\geq \frac{1}{24} \frac{(n+1)}{ (n-1)^2}$, for $n\geq 3$.

 \begin{pf} of Theorem \ref{N3.1}

Now we are ready to prove the estimate (\ref{3.3}).  We assume $\partial \hat{B}(x_0,a)$ is not empty. 
Let $h(x)=|\nabla \log u|^2(x)-\kappa \frac{n-2}{n-1}u^{-2}|E|^2+\kappa'\frac{n+1}{(n-1)^2}|d\phi|^2$, $f(x)=h(x)d_{\hat{g}}^2(x,\partial \hat{B}(x_0,a))$, and $f(\bar{x})=\sup_{x\in \hat{B}(x_0,a)}f(x)$ for some $\bar{x}\in \hat{B}(x_0,a)$.

 Note that the function $f(x)$ is invariant under the scaling of the metric.    Scaling $u$ and the metric $\bar{g}$ by $u(\bar{x})^{-1}$ and $h(\bar{x})$,   we may assume  $u(\bar{x})=1$ and $h(\bar{x})=1$.  It is not hard to prove   (see  (\ref{2.49}) (\ref{2.50})): 
 
  \begin{equation} \label{3.18}
 h(x)\leq 4 \ \ \  \text{on}  \  \hat{B}(\bar{x}, D),
  \end{equation}
    where  \begin{equation} \label{D}
       D=\frac{1}{2}\sqrt{f(\bar{x})}. 
      \end{equation}
      
From (2.16) in \cite{Chen1}, we have 
\begin{equation} \label{3.19}
\begin{split}
& \hat{R}ic({e_0,e_0})=-\bar{R}ic(e_0,e_0) \\
& \hat{R}ic({e_0,e_j})=-\bar{R}ic(e_0,e_j)=0
\\& \hat{R}ic({e_i,e_j})=\bar{R}ic(e_i,e_j).
\end{split}
\end{equation}
  Together  with (\ref{3.5}) (\ref{3.18}),  we get 
  \begin{equation} \hat{R}ic\geq -\frac{1}{n-2} h \hat{g}
  \end{equation}
  hence
  \begin{equation}\label{3.22} \hat{R}ic\geq -\frac{4}{n-2}\hat{g} \ \  \text{on}  \  \hat{B}(\bar{x}, D).
  \end{equation}
  Let $\xi:[0,\infty) \rightarrow [0,1]$ be a fixed nonnegative smooth non-increasing function such that $\xi=1$ on $[0,\frac{1}{2}]$ and $\xi=0$ on $ [1,\infty)$. Set $L(x)=\xi(\frac{d_{\hat{g}}(\bar{x},x)}{D})h(x)$, and $L(\hat{x})=\sup_{x\in  \hat{B}(\bar{x},D)}L(x)$ for some  $\hat{x}\in B(\bar{x}, D)$.  From (\ref{3.17}), we know that $h$ satisfies 
\begin{equation} \label{3.23}
\hat{\triangle} h\geq \frac{1}{72}h^2.
\end{equation}  
  
  Using (\ref{3.22}) (\ref{3.23})  and  applying maximum principle on $L(x)$ at $\hat{x}$ as in (\ref{2.54}) (\ref{2.55}) (\ref{2.56}), we find at $\hat{x}$: 
     \begin{equation}\label{3.24}
  \begin{split}
  0\geq & \frac{1}{72}\xi  h^2-\frac{2}{D^2}\frac{(\xi^{'})^2}{\xi}h\\&+h[\frac{1}{D^2}(\xi^{''}+2n\xi^{'})+\frac{2}{D}\xi^{'}\sqrt{\frac{n}{n-2}}].
  \end{split}
  \end{equation}
  Multiplying  both sides of (\ref{3.24}) by $\xi(\hat{x})$,  we find 
  $$
  \frac{1}{72}L(\hat{x})-\frac{Cn}{D^2}-\frac{C}{D}\leq 0,
  $$
  which implies 
  \begin{equation}\label{3.25}
  1\leq C(\frac{n}{D^2}+\frac{1}{D}),
  \end{equation}
   since $L(\hat{x})\geq L(\bar{x})=1$. Hence   $D\leq C\sqrt{n}$ for some universal constant $C$. From  (\ref{D}), we get the desired estimate $f(\bar{x})\leq Cn$.
   
   If $\partial \hat{B}(x_0,a)$ is empty, we can  apply the strong maximum principle on equation (\ref{3.17}) directly to conclude   $h\equiv 0$. 

   \end{pf}
   
  A corollary of Theorem \ref{N3.1} is the following theorem, which is a generalization of  Theorem 1.2 in  \cite{Chen1}.  

    \begin{thm} \label{N3.2}  Let $(M,g_{M})$ be a static spacetime of dimension $n+1$ satisfying the assumptions in Theorem \ref{N3.1}. If we assume $(M,g_{M})$ is geodesically complete, then the universal cover of $(M,g_{M})$ is isometric to a product $R\times N$ equipped with a product metric $-dt^2+g_{N}$, where $g_{N}$ is a complete Ricci flat Riemannian manifold. \end{thm}

Bing-Long Chen \\
 Department of Mathematics, \\
Sun Yat-sen University,\\
Guangzhou, P.R.China, 510275\\
Email: mcscbl@mail.sysu.edu.cn\\\end{document}